\documentclass[12pt]{article}
\usepackage[utf8]{inputenc}
\usepackage{amsmath, mathtools, amssymb, amsthm}
\usepackage[margin=1.5in]{geometry}
\usepackage{soul}

\usepackage{titlesec}
\titlespacing*{\section}{0pt}{0.5ex}{0.5ex}
\titlespacing*{\subsection}{0pt}{0.5ex}{0.5ex}
\titlespacing*{\subsubsection}{0pt}{0.4ex}{0.4ex}

\makeatletter
\renewcommand*\env@matrix[1][\arraystretch]{%
  \edef\arraystretch{#1}%
  \hskip -\arraycolsep
  \let\@ifnextchar\new@ifnextchar
  \array{*\c@MaxMatrixCols c}}
\makeatother

\title{On entropy stable temporal fluxes}
\author{Ayoub Gouasmi, University of Michigan, \\
 Karthik Duraisamy, University of Michigan, and \\
 Scott Murman, NASA Ames Research Center.}
\date{}

\begin{document}

\maketitle

\begin{abstract}
    Entropy-stable (ES) schemes have gained considerable attention over the last decade, especially in the context of turbulent flow simulations using high-order methods. While promising because of their nonlinear stability properties, ES schemes have to address a number of issues to become practical. One of them is how much entropy should be produced by the scheme at a certain level of under-resolution. This problem has been so far studied by considering different ES interfaces fluxes in the spatial discretization only because they can be tuned to generate a certain amount of entropy. In this note, we point out that, in the context of space-time discretizations, the same applies to ES interface fluxes in the temporal direction.  
\end{abstract}

\section{Motivations}
\hspace{0.45 cm} A fundamental challenge in the numerical simulation of high Reynolds number turbulent flows is that most of the  discontinuities encountered by the numerical scheme are not physical but a product of under-resolution. Classical CFD algorithms such as the Roe scheme \cite{Roe} are designed to handle ``physical discontinuities" such as shocks. These schemes do not provide any stability guarantee in the general case. Furthermore, in the context of Large Eddy Simulation (LES) they can render sub-grid-scale (SGS) models inactive \cite{Garnier}, and thus complicate the formulation and solution. \\ 
\indent The family of Entropy-Conservative (EC) and Entropy-Stable (ES) schemes introduced by Tadmor \cite{Tadmor} has gained significant  attention over the last decade. EC/ES schemes were initially developed for hyperbolic systems of conservation laws with a generalized entropy function (referred to as ``entropy" for simplicity). The Euler equations are one such example, with $\rho S$ as one admissible entropy \cite{Harten}. EC schemes solve the original system in such a way that an additional conservation equation for an entropy variable is also satisfied. An ES scheme adds dissipation to an EC scheme to enforce entropy production across interfaces with jumps. Globally, this results in a stability statement that resembles the second principle of thermodynamics. This stability property carries over to the Navier--Stokes equations provided that the additional terms (viscous diffusion and heat conduction, which physically add entropy to the system) are discretized in a consistent manner \cite{Zhong}. Several developments followed among which extensions to high-order finite-volume \cite{LeFloch, Fjordholm}, finite-difference \cite{Fisher} and finite-element \cite{Hughes, Barth} methods. Discontinuous-Galerkin (DG) methods are of particular interest in this note. Entropy Stability in DG methods can be obtained either by discretizing the so-called entropy variables \cite{Tadmor} instead of the conservative variables \cite{Barth} or by using summation-by-parts operators \cite{FisherCarpenter}. In both approaches, an entropy-stable interface flux is required. An ES space-time DG code, called \textit{eddy}, based on the former approach has been developed by Murman and co-workers \cite{Murman1, Murman2, Murman3, Murman4} over the last few years to perform scale-resolving simulations of separated flows.  \\
\indent ES schemes have yet to overcome a number of hurdles before successfully being put into practice. The second principle of thermodynamics does not specify \textit{by how much} the entropy of a closed system should increase. In the same vein, it has yet to be established how much entropy an ES scheme should produce when dealing with discontinuities, whether arising from physics or under-resolution. Too much entropy will lead to damped results, whereas too little will lead to oscillatory solutions. Diosady and Murman \cite{Murman3} compared several entropy stable fluxes in under-resolved simulations of the Taylor Green Vortex at various Reynolds number. Ismail and Roe \cite{IsmailRoe} introduced the concept of "entropy consistent" schemes for shocks in the Burgers and Euler equations. Their work focused on properly tuning the dissipation matrix of the convective flux to achieve the right entropy production across shocks. Derigs \textit{et al.} \cite{Derigs} further examined the design of the entropy producing dissipation matrix in the context of astrophysical simulations.  \\
\indent All of the above studies considered the effect of the spatial discretization only, and in problems where entropy is supposed to increase (because of shocks or viscous effects). In a recent preprint \cite{Gouasmi}, the authors examined the behavior of EC schemes in the receding flow problem (also referred to as the \textit{123 problem} by Toro~\cite{Toro}) where two rarefactions propagate in opposite directions. Conventional finite-volume schemes are known to produce an anomalous temperature rise (termed ``overheating") at the symmetry point. In a series of papers, Liou \cite{Liou1, Liou2, Liou3} established a link between the overheating and a spurious entropy rise initiated at the first time iteration. This analysis also suggested that conserving entropy would prevent the overheating. In \cite{Gouasmi}, numerical experiments using an EC flux in space highlighted the non-negligible impact of time-integration \cite{TadmorAN}. Further investigation lead the authors to develop a new time integration scheme (as part of an unsuccessful attempt at curing the overheating) that introduces no production of entropy. This scheme was derived by leveraging the analogy between the condition that must be satisfied by an EC flux in space \cite{Tadmor} and the condition which defines the EC time scheme of Lefloch \textit{et. al} \cite{LeFloch}. \\
\indent In space-time DG schemes, an upwind flux in time is typically used because of causality and also because it allows for the solution of one-time slab at a time. Upwinding in time is entropy stable \cite{Barth}, meaning that it does contribute to the damping of information. We show, using the 1D Euler equations as an example, how EC and ES temporal fluxes can be derived following the methodology used for spatial fluxes.

\section{Preliminaries}
Consider a hyperbolic system of conservation laws: 
\begin{equation}\label{eq:base}
\frac{\partial u}{\partial t} + \frac{\partial f}{\partial x} = 0, \ u = u(x,t) \in \mathbb{R}^{\mathbb{N}}, \ x \in \mathbb{R}, t > 0,
\end{equation}
which admits a convex extension \cite{Lax} in the sense that an additional conservation law: 
\begin{equation}\label{eq:base_entropy}
\frac{\partial U}{\partial t} + \frac{\partial F}{\partial x} = 0,
\end{equation}
where $U = U(u) \in \mathbb{R}$ is a convex function of $u$ and $F = F(u)$, naturally results from eq. (\ref{eq:base}). Define the vector of entropy variables $v = (\frac{\partial U}{\partial u})^T$. As a reference, we consider the one-dimensional Euler equations:
\begin{equation}\label{eq:Euler}
    \frac{\partial}{\partial t} \begin{bmatrix}
    \rho \\ \rho u \\ \rho e^t
    \end{bmatrix} +
    \frac{\partial}{\partial x}
    \begin{bmatrix}
    \rho u \\ \rho u^2 + p \\ \rho u h^t
    \end{bmatrix} = 0,
\end{equation}
where $\rho$ is the density, $u$ is the velocity, $e^t= e + \frac{1}{2}u^2$ is the total energy per unit mass, $e$ is the internal energy, $\ p = \rho (\gamma - 1) e$ is the thermodynamic pressure, $h^t = e^t + \frac{p}{\rho} = \frac{a^2}{\gamma - 1} + \frac{1}{2}u^2$ is the total enthalpy and $a = \sqrt{\frac{\gamma p}{\rho}}$ is the speed of sound. An additional conservation equation that results from eq. (\ref{eq:Euler}) is that of entropy $\rho S$ ($S$ is the specific entropy):
\begin{equation*}
    \frac{\partial (\rho S)}{\partial t} + \frac{\partial (\rho u S)}{\partial x} = 0.
\end{equation*}
For shocks, the conservation equation for entropy is replaced with:
\begin{equation*}
    \frac{\partial (\rho S)}{\partial t} + \frac{\partial (\rho u S)}{\partial x} < 0.
\end{equation*}
which should hold in the sense of distributions. We employ the entropy-entropy flux pair introduced by Hughes \textit{et al.} \cite{Hughes}:
\begin{equation*}
    U(u) = -\frac{\rho S}{\gamma - 1}, \ F(u) = -\frac{\rho u S}{\gamma - 1}, \ S = ln(p) - \gamma ln(\rho). 
\end{equation*}
The corresponding entropy variables are given by:
\begin{equation}\label{eq:EntropyVar}
v = \bigg[\frac{\gamma - S}{\gamma - 1} - \frac{1}{2}\frac{\rho u^2}{p}, \ \frac{\rho u}{p}, \ -\frac{\rho}{p}\bigg]^T.
\end{equation}
EC schemes \cite{Tadmor} are finite-volume discretizations of eq. (\ref{eq:base}):
\begin{equation}\label{eq:FVM}
\frac{d}{dt}u_{j}(t) + \frac{1}{\Delta x}[f_{j + \frac{1}{2}} - f_{j - \frac{1}{2}}] = 0,
\end{equation}
which imply a finite-volume discretization of eq. (\ref{eq:base_entropy}):
\begin{equation}
\frac{d}{dt}U(u_{j}) + \frac{1}{\Delta x}[F_{j + \frac{1}{2}} - F_{j - \frac{1}{2}}] = 0.
\end{equation}
Tadmor \cite{Tadmor} showed that an interface flux $f^*$ is EC if and only if across the interface:
\begin{equation}\label{eq:ECflux}
    [v] \cdot f^* = [\psi], \ \psi = v \cdot f - F
\end{equation}
$\psi = \rho u$ in the Euler equations with our choice of entropy variables. The first EC flux was proposed by Tadmor \cite{Tadmor}:
\begin{equation}\label{eq:TadmorEC}
f_{j + \frac{1}{2}} = \int_{0}^1 f(v_{j + \frac{1}{2}}(\xi)) d\xi, \  v_{j + \frac{1}{2}}(\xi) = v_{j} + \xi \Delta v_{j + \frac{1}{2}}, \ \Delta v_{j + \frac{1}{2}} = v_{j + 1} - v_{j}.
\end{equation}
This flux has the inconvenient property of not having a closed form. Subsequently, Roe \cite{RoeEC} proposed a simpler flux $f^* = [f_1, \ f_2, \ f_3]$ for the Euler equations:
\begin{equation}\label{eq:RoeEC}
f_1 = \bar{z_2} z_3^{ln} , \
f_2 = (\bar{z_3} + f_1 \bar{z_2})/(\bar{z_1}) , \
f_3 = \frac{1}{2 \bar{z_1}}(- f_1 \frac{1 + \gamma}{1 - \gamma}\frac{1}{z_1^{ln}} + f_2 \bar{z_2}),
\end{equation}
where $z_1 = \sqrt{\frac{\rho}{p}}$, $z_2 = \sqrt{\frac{\rho}{p}}u$ and $z_3 = \sqrt{\rho p}$ are independent variables used to algebraically solve for $f^*$ in eq. (\ref{eq:ECflux}) \cite{RoeEC, IsmailRoe}. $\bar{z}$ and $z^{ln}$ denote the arithmetic and logarithmic averages, respectively. \\
\indent ES schemes are finite-volume discretizations of eq. (\ref{eq:base}) which imply:
\begin{equation}\label{eq:ES}
\frac{d}{dt}U(u_{j}) + \frac{1}{\Delta x}[F_{j + \frac{1}{2}} - F_{j - \frac{1}{2}}] = -\mathcal{E} < 0.
\end{equation}
Tadmor \cite{Tadmor} showed that by combining an EC flux with a dissipation term in the interface flux:
\begin{equation}
f_{j+\frac{1}{2}} = f_{j + \frac{1}{2}}^* - Q_{j+\frac{1}{2}}\Delta v_{j+\frac{1}{2}},
\end{equation}
with $Q_{j+\frac{1}{2}}$ positive definite, equation (\ref{eq:ES}) is enforced with:
\begin{equation}\label{eq:SpaceDiss}
\mathcal{E} = \frac{1}{2 \Delta x}\big[ \Delta v_{j + \frac{1}{2}}^T Q_{j + \frac{1}{2}} \Delta v_{j + \frac{1}{2}} + \Delta v_{j - \frac{1}{2}}^T Q_{j - \frac{1}{2}} \Delta v_{j - \frac{1}{2}} \big] > 0.
\end{equation}
Ismail and Roe \cite{IsmailRoe} developed a Roe-type dissipation matrix using Barth's eigenscaling theorem \cite{Barth}.
\section{Entropy stable temporal fluxes}
In a space-time DG method where the entropy variables are discretized and quadrature errors are assumed to be negligible, entropy stability depends solely on the interface fluxes used in space and time. Hence it suffices to consider a space-time finite volume discretization of eq. (\ref{eq:base}):
\begin{equation}\label{eq:STFVM}
[u_{j}^{n+\frac{1}{2}} - u_{j}^{n - \frac{1}{2}}] + \lambda [f_{j + \frac{1}{2}}^n - f_{j - \frac{1}{2}}^n] = 0, \ \lambda = \frac{\Delta t}{\Delta x}
\end{equation}
In cell $(j, n)$, $u_{j}^{n}$ is the mean solution. $f_{j + \frac{1}{2}}^{n}$ is a consistent spatial flux. $u_{j}^{n + \frac{1}{2}}$ is a consistent temporal flux. Entropy conservation is achieved if eq. (\ref{eq:STFVM}) implies a similar discretization for the entropy equation (\ref{eq:base_entropy}):
\begin{equation}
[U_{j}^{n+\frac{1}{2}} - U_{j}^{n - \frac{1}{2}}] + \lambda [F_{j + \frac{1}{2}}^n - F_{j - \frac{1}{2}}^n] = 0
\end{equation}
with $F_{j + \frac{1}{2}}^n$ a consistent entropy spatial flux and $U_{j}^{n + \frac{1}{2}}$ a consistent entropy temporal flux. Likewise, the scheme is entropy-stable if:
\begin{equation}
[U_{j}^{n+\frac{1}{2}} - U_{j}^{n - \frac{1}{2}}] + \lambda [F_{j + \frac{1}{2}}^n - F_{j - \frac{1}{2}}^n] < 0
\end{equation}
is inherently solved. 
Entropy conservation in time is achieved if: 
\begin{equation}
(v_{j}^n)^T [u_{j}^{n + \frac{1}{2}} - u_{j}^{n - \frac{1}{2}}] = U_{j}^{n + \frac{1}{2}} - U_{j}^{n - \frac{1}{2}}.
\end{equation}
From this point, we can re-use the analysis of Tadmor \cite{Tadmor} and show that an EC flux in time $u^*$ should satisfy:
\begin{equation}\label{eq:ECTcond1}
[v] \cdot u^* = [\phi], \ \phi = v \cdot u - U,
\end{equation}
where $\phi$ is the temporal flux potential. For the Euler equations and the specific entropy-entropy flux pair in question, $\phi = \rho$. From there, we can derive EC temporal fluxes that are the time counterparts of EC spatial fluxes. The temporal version of Tadmor's EC temporal flux \cite{Tadmor} is given by:
\begin{equation}\label{eq:TadmorTEC}
u_{j}^{n + \frac{1}{2}} = \int_0^1 u(v_{j}^n + \xi \Delta v_{j}^{n + \frac{1}{2}}) d\xi, \ \Delta v_{j}^{n + \frac{1}{2}} = v_{j}^{n+1} - v_{j}^{n}.
\end{equation}
To derive the counterpart of Roe's EC flux in time, we adopt the same approach \cite{RoeEC, IsmailRoe}.  We can use the same algebraic variables $z_i$ as for the spatial flux. Denoting $u^* = [u_1, \ u_2, \ u_3]$ , eq. (\ref{eq:ECTcond1}) can be rewritten as:
\begin{equation*}
    u_1 ( \frac{1}{z_3^{ln}} [z_3] - \frac{1 + \gamma}{1 - \gamma}\frac{1}{z_1^{ln}} [z_1] - \bar{z_2}[z_2]) +
    u_2 (\bar{z_1}[z_2] + \bar{z_2}[z_1]) + u_3 (- 2 \bar{z_1}[z_1]) = \bar{z_1}[z_3] + \bar{z_1}[z_3].
\end{equation*}
Regrouping, we get:
\begin{equation*}
    [z_1](- u_1 \frac{1 + \gamma}{1 - \gamma}\frac{1}{z_1^{ln}} + u_2 \bar{z_2} - 2 u_3 \bar{z_1}) + [z_2] ( - u_1\bar{z_2} + u_2 \bar{z_1}) + [z_3](\frac{1}{z_3^{ln}} u_1)
     = [z_1] \bar{z_3} + [z_3]\bar{z_1}.
\end{equation*}
The jumps in the $z_i$ are independent, therefore:
\begin{equation*}
- u_1 \frac{1 + \gamma}{1 - \gamma}\frac{1}{z_1^{ln}} + u_2 \bar{z_2} - 2 u_3 \bar{z_1} = \bar{z_3}, - u_1\bar{z_2} + u_2 \bar{z_1} = 0, \ \frac{1}{z_3^{ln}} u_1 = \bar{z_1}. 
\end{equation*}
The resulting temporal flux is given by:
\begin{equation}\label{eq:RoeECT}
u_1 = \bar{z_1} z_3^{ln}, \
u_2 = u_1 \frac{\bar{z_2}}{\bar{z_1}}, \
u_3 = \frac{1}{2 \bar{z_1}}(- u_1 \frac{1 + \gamma}{1 - \gamma}\frac{1}{z_1^{ln}} + u_2 \bar{z_2} - \bar{z_3}).
\end{equation}
The extension to higher dimensions and other hyperbolic systems with a convex extension is straightforward. \\
\indent Entropy stability in time can be achieved in the same manner as in the spatial case, namely by adding a dissipation term to an entropy conservative temporal flux. It can be shown, by again following Tadmor \cite{Tadmor} that a dissipation term of the form $T^{n + \frac{1}{2}} \Delta v^{n + \frac{1}{2}}$, with $T^{n + \frac{1}{2}}$ a positive definite matrix, qualifies, and that the entropy production in cell $n$ due to temporal fluxes will be given by:
\begin{equation}\label{eq:TimeDiss}
\mathcal{E} = \frac{1}{2 \Delta t}\big[ (\Delta v^{n + \frac{1}{2}})^T T^{n + \frac{1}{2}} \Delta v^{n + \frac{1}{2}} + (\Delta v^{n - \frac{1}{2}})^T T^{n - \frac{1}{2}} \Delta v^{n - \frac{1}{2}} \big].
\end{equation}
The jacobian $H(v)$ of the temporal flux $u(v)$ with respect to the entropy variables $v$ is positive definite, it can therefore be used as a $T^{n+\frac{1}{2}}$. 

\noindent We conclude this note by providing an alternative proof that upwinding in time is entropy stable. The proof relies on rewriting upwinding in time as the combination of an EC flux and a dissipation term. Upwinding in time is given by $u^{n + \frac{1}{2}} = u^{n}$. An argument for using this flux is causality, which is a result of  the positive definiteness of the temporal jacobian. The upwind flux can be re-written as: 
\begin{equation*}
    u^{n} = \frac{u^{n} + u^{n+1}}{2} - \frac{1}{2} (u^{n+1} - u^{n}).
\end{equation*}
The jump term can be rewritten as:
\begin{equation}\label{eq:jump_u}
    u^{n+1} - u^{n} = \int_{v^{n}}^{v^{n+1}} H(v) dv.
\end{equation}
The integral in eq. (\ref{eq:jump_u}) is independent of the path chosen. Taking a straight line $v(\xi) = v^{n + \frac{1}{2}}(\xi) = v^{n} + \xi \Delta v^{n + \frac{1}{2}} $ gives:
\begin{equation}\label{eq:jump_u2}
    u^{n+1} - u^{n} = \int_{0}^{1} H(v^{n + \frac{1}{2}}(\xi)) \Delta v^{n+\frac{1}{2}} d\xi.
\end{equation}
Equation (\ref{eq:jump_u2}) allows us to rewrite the upwind flux in ``viscosity form" \cite{Tadmor}:
\begin{equation}\label{eq:upwindTvisc}
    u^{n} = \frac{u^{n} + u^{n+1}}{2} - \bigg [\int_{0}^{1} \frac{1}{2} H(v^{n + \frac{1}{2}}(\xi))  d\xi \bigg ] \Delta v^{n+\frac{1}{2}}
\end{equation}
The temporal version of Tadmor's EC flux  (\ref{eq:TadmorTEC}), denoted $u^*$ in this context, also has a "viscosity form" that is analogous to that of the spatial version (\ref{eq:TadmorEC}):
\begin{equation}\label{eq:TadmorECTvisc}
u^* = \frac{u^{n} + u^{n+1}}{2} - \bigg [ \int_0^1 (\xi - \frac{1}{2}) H(v^{n + \frac{1}{2}}(\xi)) d\xi\bigg] \Delta v^{n + \frac{1}{2}}.
\end{equation}
Comparing the viscosity forms (\ref{eq:upwindTvisc}) and (\ref{eq:TadmorECTvisc}), we can finally rewrite upwinding in time as:
\begin{equation}\label{eq:UpwindDecomp}
    u^n = u^* - T^{n + \frac{1}{2}} \Delta v^{n + \frac{1}{2}}
\end{equation}
where $T^{n + \frac{1}{2}} = \bigg [ \int_0^1 (1 - \xi) H(v^{n + \frac{1}{2}}(\xi)) d\xi\bigg]$ is symmetric positive definite. This completes the proof. 

\noindent {\em Note 1:} Since the temporal jacobian is positive definite, upwinding in time will \textit{always} produce entropy in the presence of a jump in states. 

\noindent {\em Note 2:} The dissipation term in the decomposition (\ref{eq:UpwindDecomp}) of upwinding in time can be seen as a ``temporal" implicit SGS model \cite{Boris1, Boris2}.  
\section{Additional remarks}
\hspace{0.5 cm} In a space-time DG formulation, the finite element solution in the time slab $t^n$ is determined by the fluxes across the $n+\frac{1}{2}$ and $n-\frac{1}{2}$ temporal interfaces, in addition to the fluxes in the spatial directions (note that, in the case of a moving geometry \cite{Murman4}, which space-time DG methods can naturally adjust to, the spatial and temporal interfaces are no longer orthogonal to each other). With upwinding in time, the solution in time slab $t^n$ is only influenced by the solution in time slab $t^{n-1}$. The EC flux (\ref{eq:RoeECT}) we derived in this paper is clearly non-causal and would couple all the time slabs together. The complete coupling between time slabs can be avoided if the non-causal temporal flux is used in the $n - \frac{1}{2}$ interface only. One could also solve one block of consecutive, coupled time slabs over a given time interval by appropriately choosing the fluxes (causal or non-causal) at the interfaces. From a computational viewpoint, this approach incurs higher memory costs but also a higher arithmetic intensity which could make space-time DG schemes even more amenable to parallel computing. \\   
\indent The goal of this note is not to question the validity of causality as a physical principle, but rather to ask whether relaxing causality with entropy stability could bring about some improvements in the development of accurate ES space-time schemes for turbulent flows. If upwinding in time appears to be indispensable, then another key question to answer is: To what extent is entropy production in space required, considering that entropy is already being produced in time?

\section*{Acknowledgments}
This work was partially supported by the AFOSR through grant no. FA9550
-16-1-030.

\end{document}